\documentclass[11pt]{amsart}

\usepackage{amsfonts,amsmath,amscd,latexsym}
\input xy
\xyoption{all}

\def\PP{{\textbf P}}
\def\OO{{\mathcal O}}
\def\II{{\mathcal I}}

\def\H{\mathcal{H}}

\def\Pic0{{\rm Pic}^0(X)}

\theoremstyle{plain}

\newtheorem{theoremalpha}{Theorem}

\newtheorem{theorem}{Theorem}[section]
\newtheorem{proposition}[theorem]{Proposition}
\newtheorem{corollary}[theorem]{Corollary}
\newtheorem{lemma}[theorem]{Lemma}
\newtheorem{lemma/notation}[theorem]{Lemma/Notation}

\theoremstyle{definition}
\newtheorem{definition}[theorem]{Definition}
\newtheorem{remark}[theorem]{Remark}
\newtheorem{example}[theorem]{Example}

\newtheorem{conjecture}[theorem]{Conjecture}
\newtheorem{conjecture/question}[theorem]{Conjecture/Question}

\newtheorem{remark/definition}[theorem]{Remark/Definition}
\newtheorem{notation}[theorem]{Notation}
\newtheorem{problem}[theorem]{Problem}
\newtheorem{definition/notation}[theorem]{Definition/Notation}

 \theoremstyle{remark}

\begin{document}

\title{Castelnuovo theory and the geometric Schottky problem}

\author[G. Pareschi]{Giuseppe Pareschi}
\address{Dipartimento di Matematica, Universit\`a di Roma, Tor
Vergata, V.le della Ricerca Scientifica, I-00133 Roma, Italy}
\email{{\tt pareschi@mat.uniroma2.it}}

\author[M. Popa]{Mihnea Popa}
\address{Department of Mathematics, University of Chicago,
5734 S. University Ave., Chicago, IL 60637, USA } \email{{\tt
mpopa@math.uchicago.edu}}

\thanks{The second author was partially supported by the NSF grants DMS-0200150
and DMS-0500985, and by an AMS Centennial Fellowship.}

\date{\today}
\maketitle

\setlength{\parskip}{.1 in}

\section{Introduction}

The aim of this paper is to show that Castelnuovo theory in projective space (cf. \cite{acgh} Ch.III \S2 and \cite{gh} Ch.4
\S3) has a precise analogue for abelian varieties. This can be quite surprisingly related in a concrete way to the geometric Schottky problem, namely the problem of identifying Jacobians among all principally polarized abelian varieties (ppav's) via geometric conditions on the polarization. The main result is that a ppav satisfies a direct analogue of the Castelnuovo Lemma if and only if it is a Jacobian. We prove or conjecture other results which show an extremely close parallel between geometry in projective space and Schottky-type projective geometry on abelian varieties.

On a ppav $(A, \Theta)$ of dimension $g$ one can make sense of what it means for a finite set $\Gamma$ of at least $g+1$ distinct
points to be in general position (we call this \emph{theta-general position}, cf. \S3): we simply require for any subset
$Y\subset \Gamma$ of $g+1$ points to be \emph{theta-independent}, i.e. for any  $g$ points of $Y$ there is a translate of $\Theta$
containing them and avoiding the remaining point. It turns out that general points on any non-degenerate curve are in
theta-general position.
 On a Jacobian $J(C)$, points on an Abel-Jacobi curve curve $C$  impose the minimal
number of conditions, namely $g+1$, on the linear series $|(2\Theta)_\alpha|$, for $\alpha\in J(C)$ general (Proposition
\ref{conditions} and Example \ref{ajexample}). The main result we prove is the following theorem -- we refer to \S2 for a detailed description of the context in which it should be integrated.

\begin{theoremalpha}[``Castelnuovo-Schottky Lemma"]\label{alpha}
Let $(A,\Theta)$ be an irreducible principally polarized abelian variety of dimension $g$, and let $\Gamma$ be a set of
$n\geq g+2$ points on $A$ in theta general position, imposing only $g+1$ conditions on the linear series
$|\OO_A(2\Theta)\otimes\alpha|$ for $\alpha$ general in $ \widehat A$. Then $(A, \Theta)$ is the canonically polarized
Jacobian of a curve $C$ and $\Gamma \subset C$ for a unique Abel-Jacobi embedding $C\subset J(C)$.
\end{theoremalpha}

Roughly speaking, the key points in the proof of the Theorem are the following. Given a set  of points $\Gamma$ in theta-general position and imposing generically only $g+1$ conditions on $|\OO_A(2\Theta)\otimes\alpha|$, and given a subset
$Z\subset \Gamma$ with $|Z|=g+1$, we consider the locus $V(Z)$ of $\alpha \in \widehat A$  such that $Z$ fails to impose independent conditions on $|\OO_A(2\Theta)\otimes\alpha)|$.  It turns out that  $V(Z)$ is a \emph{theta-translate}, precisely described in function of $Z$ (cf. Proposition \ref{keyresult}). We prove a formula for the intersection of theta-translates of type $V(Z)$ (cf. Proposition \ref{elementary}). The intersection formula yields the existence of a certain positive-dimensional family of trisecants to the Kummer variety (cf. Theorem \ref{trisecants}), at which stage the Gunning-Welters criterion \cite{welters1} implies that $(A,\Theta)$ is a Jacobian. This approach also carries a natural way to recover the curve.
In fact, in analogy with the classical Castelnuovo setting, the  curve $C$ is
 recovered as the base locus of a continuous system of divisors algebraically equivalent to $2\Theta$,
 containing $\Gamma$ (cf. Corollary \ref{torelli}). In particular this provides another proof of Torelli. An
 important ingredient is a previous result on the regularity of an Abel-Jacobi curve (cf. \cite{pp1}, \S4). 

In keeping with the picture suggested by Castelnuovo theory, we also establish a genus bound for non-degenerate curves in ppav's.

\begin{theoremalpha}\label{genus_bound}
Let $(A,\Theta)$ be a $g$-dimensional irreducible ppav. Let $C^\prime$ be a smooth curve of genus $\gamma$, admitting a birational map onto a non-degenerate curve $C$ of degree $d: = C\cdot \Theta$ in $A$. Let
$m=\bigl[\frac{d-1}{g}\bigr]$, so that $d-1=mg+\epsilon$,  with $0\le \epsilon<g$. Then
$$\gamma\le {{m+1}\choose 2}g+(m+1)\epsilon+1.$$
Moreover, the inequality is strict for $g\ge 3$ and $d\ge g+2$.
\end{theoremalpha}

The bound is quadratic in the degree, of leading term $\frac{d^2}{2g}$, so of the order suggested by curves in projective space. It is optimal for
abelian surfaces, but the Castelnuovo-Schottky Lemma implies that in higher dimensions
-- unlike in the case of projective space -- there is room for improvement (cf. \S6).

We defer for the next section a detailed discussion of the context described in this Introduction, as well as of conjectural developments and connections with our previous work concerning regularity. Let us mention here that, as a consequence of criteria involving $M$-regularity, in \S7 we show how  to attach a canonical \emph{divisor class} to a uniform collection of points $\Gamma$ failing to impose independent conditions on $2$-theta functions generically. This concept seems to hold the key for future developments.

Finally, we mention that the analogue of Theorem \ref{alpha} in the case of finite schemes, in the spirit of the Eisenbud-Harris generalization
of the Castelnuovo Lemma (\cite{eh2}, \cite{eh3}), will appear in work of M. Lahoz \cite{lahoz}.

\noindent \textbf{Acknowledgements.} We would like to thank Ciro
Ciliberto, Olivier Debarre, Lawrence Ein, Joe Harris and Elham
Izadi, for their interest, suggestions and answers to numerous
questions. We are grateful to Mart\'i Lahoz and Juan Carlos Naranjo for pointing out an error
in a previous version of the paper. Special thanks are due to Sam
Grushevsky for correspondence and for communicating to us his work
\cite{grushevsky1}, which contains a statement close to our first
theorem (cf. \S5.2 below for the precise statement and comparison).
We would also like to emphasize that we realized our main theorem
might also work for $g+2$ points, as opposed to $g+3$ or more, only
after seeing his paper.

 \section{A parallel between projective spaces and principally polarized abelian varieties}

The starting point of our work is the observation that on ppav's
there is an extensive similarity between basic facts related to the
geometry of points in general linear position and of rational normal
curves in projective space on one hand, and the geometry of points
in theta-general position and of Abel-Jacobi embedded curves in
Jacobians, on the other hand.  This begins with a \emph{point-theta
divisor} correspondence, similar to the point-hyperplane
correspondence between a
 projective space and its dual, and continues as described shortly. Ideally one then hopes that on ppav's there are
 analogues of those aspects of projective geometry which are consequences of the ``geometry of hyperplanes". Among these
 aspects there is Castelnuovo theory.\footnote{Castelnuovo theory initiated with Castelnuovo's work \cite{castelnuovo}
 around 1890. Modern accounts, as well as new results and research directions, have been given, among others, in
  \cite{gh}, \cite{eh1}, \cite{acgh}, \cite{ciliberto}, \cite{reid}, \cite{egh}. The current perspective on the subject is mostly due to J. Harris.} It becomes natural to expect an analogue of Castelnuovo theory for ppav's, at least in its basic aspects. What is perhaps more surprising is that this turns out to lead to a geometric characterization of Jacobians (Theorem \ref{alpha} in the Introduction) among all ppav's.

Another -- this time homological -- fact pointing towards such a parallel stems from previous work in which we have explored an abelian varieties analogue of Castelnuovo-Mumford regularity (cf. \cite{pp4} for an overview of this circle of ideas). This yields results on the geometry of abelian varieties and their subvarieties which parallel
classical facts of projective geometry (\cite{pp1},\cite{pp2}, \cite{pp3}, \cite{debarre3}). Although in this paper we will mainly use elementary geometric methods and Jacobian criteria based on the existence of trisecants to the Kummer, it is our hope that they will eventually naturally combine with homological methods. A first step in this direction is made in \S8.

We list below a few entries in this analogy, taking also the
opportunity to introduce some notation. We recall that a ppav
$(A,\Theta)$ is said to be \emph{irreducible} if the theta-divisor
$\Theta$ is irreducible (as it is well known,  this means that
$(A,\Theta)$ is not isomorphic, as polarized variety, to the product
of lower dimensional ppav's). Let $(A, \Theta)$ be
 an irreducible ppav of dimension $g$, and, without loss of generality, let us assume that $\Theta$ is symmetric. For a smooth
  projective curve $C$ of genus $g$, let $(J(C), \Theta)$ be its Jacobian with the canonical principal polarization.

\noindent \textbf{(1).[Point-divisor correspondence]} (a) On $\PP^g$: the family of hyperplanes is parametrized by another
projective space of the same dimension, the dual projective space ${\PP^g}^*$. Points of $\PP^g$ correspond to hyperplanes
 of ${\PP^g}^*$, via $p\mapsto D(p)=\{[H]\in {\PP^g}^*\>|\>p\in H\}$ and $D\mapsto p(D)=\bigcap_{[H]\in D}H$.\\
 (b) On $A$: the family of divisors algebraically equivalent to $\Theta$ is parametrized by the dual variety $\widehat A$:
 \ $\{\Theta_\alpha\}_{\alpha\in \widehat A}$, where $\Theta_\alpha$ denotes the unique effective divisor in
$|\OO_A(\Theta)\otimes\alpha|$. $\widehat A$ is principally polarized and there is a correspondence between points of $A$ and
theta-divisors in $\widehat A$ given by $p\mapsto W(p)=\{\alpha\in\widehat A \>|\> p\in\Theta_\alpha\}$ and $W\mapsto
p(W)=\bigcap_{\alpha\in W}\Theta_\alpha$. The divisor $\widehat \Theta:=W(0)$ is symmetric.

\begin{notation}\label{notation}
(i) The polarization $p\mapsto \OO_A(\Theta_{p}-\Theta)$ provides the identification
 $\Psi: (A,\Theta)\rightarrow (\widehat
A,\widehat\Theta)$ so the $\Theta_\alpha$'s are translates of
$\Theta$. The translate $\Theta_{p}$ (on $A$) is identified to the
divisor ${\widehat\Theta}_p$ (on $\widehat A$), where $p$ is
identified to a line bundle on $\widehat A$. We will denote both the
above divisors  $\Theta_p$  and we will refer to them as
\emph{theta-translates}. With such identification, given
$p\in A$ and $\alpha\in \widehat A$, we have that $\alpha\in\Theta_p$ if and only if $p\in\Theta_\alpha$.  \\
(ii) Given a subscheme $Y\subset A$, we denote more generally
$$W(Y) :=\{\alpha\in \widehat A\>|\>
h^0(A,\II_Y(\Theta)\otimes\alpha)>0\>\}.$$ This locus parametrizes the theta-translates containing $Y$.\footnote{ It can be
shown that $W(Y)$ is equipped with a natural scheme structure.}  As $\Theta$ is assumed to be symmetric, if $Y=\{p\}$ then
$W(p)$ is identified to the theta-translate
$\Theta_p$ (cf. (i)).\\
(iii) Moreover we will denote
$$\OO((k\Theta)_{\alpha}) := \OO(k\Theta)\otimes \alpha.$$
\end{notation}

\noindent \textbf{(2).[General position and bound on the number of linear conditions]}
(a) On $\PP^g$: by Castelnuovo's basic remark, based on reduction to the variety of  $k$-forms which are product of linear
ones, any linearly general subset $\Gamma$ of $\PP^g$ imposes at least ${\rm min}\{|\Gamma|,kg+1\}$ independent conditions on forms of degree $k$ (cf. e.g. \cite{gh} p.252).\\
(b) On $A$: one can make sense of the notion of linear generality (\emph{theta-generality},
Definition \ref{theta-general} below). From an argument close to Castelnuovo's it follows that any theta-general subset $\Gamma$
imposes at least ${\rm min}\{|\Gamma|,(k-1)g+1\}$ conditions on $H^0(\OO_A((k\Theta)_\alpha))$ for $\alpha$ general in $\widehat A$ (cf. Proposition \ref{conditions}).

\noindent \textbf{(3).[Effectivity of the bound and curves of minimal degree]}
(a) On $\PP^g$: divisors of rational normal curves $C\subset \PP^g$ show that the bound in the previous point is sharp.
Rational normal curves are the curves of minimal degree in $\PP^g$: a \emph{non-degenerate} (i.e. not contained in a
hyperplane) curve in $\PP^g$ has degree $C\cdot H\ge g$ and
equality holds if and only if $C$ is a rational normal curve. \\
(b) On $A$: divisors of Abel-Jacobi curves show that the bound of the previous point is sharp
(Example \ref{ajexample}).  Abel-Jacobi curves are the curves of minimal degree on $A$: a curve $C$ on $A$ which is
\emph{non-degenerate}
(in the sense of groups, i.e. no translate of $C$ is contained in an abelian subvariety of $A$) has degree $C\cdot\Theta\ge g$
and equality holds if and only $C$ is an Abel-Jacobi curve. This is (a particular case of)  the Matsusaka-Ran criterion (\cite{matsusaka}, \cite{ran1}, see also \cite{bl} 11.8.1). At least partially,
it can be derived from the bound at the previous point (cf. Remark \ref{matsusaka}).

 \noindent \textbf{(4).[Castelnuovo's Lemma]}
 (a) On $\PP^g$: Castelnuovo's Lemma says that the example in (3a) is in fact the only one achieving equality:
 \emph{any set $\Gamma$ of at least $2g+3$ points in $\PP^g$ in linear general position, imposing only $2g+1$
 conditions on quadrics, lies on a unique rational normal curve.} The rational normal curve
can be recovered as the base locus of the system of quadrics through $\Gamma$
(cf. \cite{gh}, p.531).\\
(b) Under the hypothesis that $\Theta$ is irreducible, the ``Castelnuovo-Schottky Lemma"  (Theorem \ref{alpha}),
says that the example in (3b) is the only one, providing a characterization of Jacobians. The analogue of the last part of (a) supplies a Torelli-type statement:  the  Abel-Jacobi curve can be recovered as the base locus of a certain system of
divisors algebraically equivalent to $2\Theta$ passing through $\Gamma$ (cf. Corollary \ref{torelli}).

 \noindent \textbf{(5).[Castelnuovo's bound]}
(a) On $\PP^g$: Castelnuovo used the bound in (2a) to deduce his celebrated
genus bound. It turns out that the genus of curves in $\PP^g$ is bounded by a quadratic polynomial in the degree, whose leading term is
$\frac{d^2}{2(g-1)}$ (\cite{gh}, p.251--252). The argument involves the number of conditions imposed by
the general hyperplane section of a curve --  which is proved to be in linear
general position in $\PP^{g-1}$ -- to $k$-forms on $\PP^{g-1}$.\\
(b) On $A$: here  some differences arise and the results so far are not optimal.
On one hand, by Proposition \ref{thetasections} below, if the degree $d = C\cdot \Theta >g$, then a general theta-section is already theta-general. An argument similar to Castelnuovo's then shows that the genus
of a curve is bounded by a quadratic polynomial in $d$ whose leading term is $\frac{d^2}{2g}$ (Theorem \ref{genus_bound}). On the other hand, unlike in the case of projective space, the same argument together with the Castelnuovo-Schottky Lemma also shows that this bound can be improved as soon as $g\ge 3$.

\smallskip
\noindent
\textbf{Conjectural extensions.}
So far for the results of this paper, which will be addressed starting with \S3. An intriguing development would be to extend
the parallel to all \emph{varieties of minimal degree}.

\noindent \textbf{(6).[Varieties of minimal degree versus varieties representing the minimal class]}
(a) On $\PP^g$: the minimal degree of a non-degenerate subvariety of codimension $c\ge 2$ is $c+1$. Varieties of minimal
degree
are (a cone over) one of the following: (1) a rational normal scroll; (2) a Veronese surface in $\PP^5$.\\
(b) It is conjectured in \cite{debarre2} (together with some amount of evidence, including a proof on Jacobians) that on
irreducible ppav's the only subvarieties  \emph{representing the minimal class}, i.e. of codimension $d$ and of class
$[\Theta]^d/d!$, are: (1) the special varieties $W_d$ in Jacobians; (2) the Fano surface of lines in the intermediate
Jacobian of a cubic threefold.

There is a striking similarity between the two pictures (even the exceptions happen to both be surfaces in a five-fold!) and ideally there should be a geometric correspondence relating rational normal scrolls and the Veronese surface on one hand, and $W_d$'s and the Fano surface on the other hand. Finding similar properties shared by the two sets of varieties should already be important. In the next item we propose a step in this direction.

\noindent
\textbf{(7). [Cohomological regularity of the ideal sheaf]}
(a) On $\PP^g$:  a characterization of subvarieties of minimal degree is that they are the only
\emph{$2$-regular} ones. This means that the twisted ideal sheaf $\II_Y(2)$ is 0-regular, in the sense of
Castelnuovo-Mumford.\\
(b) On $A$: in \cite{pp1} Theorem 4.1 (see also \cite{pp3} Theorem
4.3) we proved that the the special subvarieties $W_d$ in Jacobians
are \emph{strongly $3$-Theta-regular}. This means that the twisted
ideal sheaf $\II_{W_d}(2\Theta))$ satisfies the Index Theorem with
index $0$, i.e. $$H^i (\II_Y((2\Theta)_\alpha)) = 0, ~\forall i>0,~
\forall \alpha\in {\rm Pic}^0(A).\footnote{We have emphasized in
\cite{pp1} that $3$-$\Theta$-regularity is the precise analogue of
Castelnuovo-Mumford $2$-regularity in projective space.}$$ The same
thing was recently checked by H\" oring \cite{horing} for the Fano surface of lines, at least for a 
general cubic threefold.  Therefore we are lead to the following:

\begin{conjecture}\label{regularity}
\emph{A non-degenerate subvariety $Y$ of an irreducible ppav $(A, \Theta)$ represents a minimal class if and only if its ideal sheaf is strongly $3$-$\Theta$-regular.}
\end{conjecture}

\subsection*{Possible extensions of the Castelnuovo-Schottky Lemma}
It is tempting to ask whether there is an interesting stratification
of the moduli space of ppav's via Castelnuovo-type conditions,
namely the existence of collections of theta-general points,
generically imposing few conditions on $2\Theta$-linear series or,
more generally, on theta-linear series of higher order. As in higher
Castelnuovo theory for projective spaces (\cite{fano}, \cite{eh1},
\cite{egh}, \cite{ciliberto},\cite{reid}) such conditions are
related, at least conjecturally, to curves of low degree. On the
other hand, there are exceptional abelian varieties containing
non-degenerate curves of very low degree (say between $g$ and $2g$)
and their geometry is quite delicate (we refer to \cite{debarre1}
for interesting results and conjectures). This suggests that for
ppav's a more careful approach is necessary, most likely based on
the \emph{existence of curves of maximal genus among those
representing a given multiple of the minimal class}. Note that such
perspective is naturally related with Prym-Tyurin theory. To this
end, we introduce in \S8 the concept of a \emph{divisor class}
attached to a uniform collection of points.

To be specific, the first natural higher Castelnuovo-Schottky problem arising is to characterize Prym varieties via a
Castelnuovo-type condition. Here the idea is suggested by a beautiful result of Welters, essentially characterizing Prym
varieties via the existence of a curve (the Abel-Prym curve) of maximal genus among those representing twice the minimal
curve class (cf. \cite{welters2}).  Let $\tilde C$ be an Abel-Prym curve in a Prym variety $(P,\Xi)$. It turns out that a
general divisor $\Gamma$ of degree $\ge 2g+1$ on $\tilde C$ is theta-general and imposes generically $n(\Gamma)=2g$
conditions on $H^0(\OO_P((2\Xi)_\alpha))$. Moreover the divisor class of $\Gamma$ is $2[\Xi]$ (see Example \ref{example}).

\begin{conjecture}\label{prym}
\emph{ Let $\Gamma$ be a theta-general,  uniform collection of points on a ppav $(A, \Theta)$, imposing generically $n<|\Gamma|$
conditions on $2\Theta$-linear series, and assume that the  divisor class associated to  $\Gamma$ is  $[2\Theta]$.
Then $n\ge 2g$ and equality characterizes Prym varieties.}
\end{conjecture}

\section{Theta general position and linear conditions on theta linear series of higher order}

\subsection{Theta general position}
In this subsection we consider natural analogues on ppav's of basic notions of linear algebra, such as linear independence and linear
general position. The main point is that, as on ppav's there is no direct analogue of linear subspaces of codimension higher than one, one is forced to define such notions using   codimension-one objects only. Let $(A, \Theta)$ be a ppav of dimension $g$.

\begin{definition/notation}\label{independence}
A collection $Z$ of $n\le g+1$ distinct points
on $A$ is \emph{theta-independent} if, for any decomposition of $Z$ as $Z=Y\cup\{p\}$, there is a theta-translate
$\Theta_{\gamma}$ such that $Y\subset \Theta_{\gamma}$ and $p\not\in \Theta_{\gamma}$. The subset of $A$ parametrizing the
family of such theta-translates is denoted $\H^{Y,p}$. The closure of $\H^{Y,p}$ is the union of some components of
$\bigcap_{q\in Y}\Theta_q$. Therefore $\dim \H^{Y,p}\ge g-n+1$ (the \emph{expected} dimension).
\end{definition/notation}

\begin{definition}\label{theta-general}  \ A collection $\Gamma$ of $n\ge g+1$  distinct points on $A$ is
\emph{theta-general} if   any $Z\subset \Gamma$, with $|Z|=g+1$, is theta-independent. In other words: for any $Y\subset
\Gamma$ with $|Y|=g$ and any $p\in \Gamma-Y$ there exists at least one theta-translate $\Theta_{\gamma}$ such that $Y\subset
\Theta_{\gamma}$ and $p\not\in \Theta_{\gamma}$.
\end{definition}

\begin{remark}\label{k-1}
All subsets of a theta-independent set are theta-independent.
Indeed, let $Z$ be theta-independent, and let $T\subset Z$. Then,
for any $q\in T$, $\H^{T,q}$ is non-empty, since obviously
$\H^{Z-\{q\},q}\subset \H^{T-\{q\},q}$.
\end{remark}

\begin{example}\label{exthetageneral}
(a) \ It is easily seen that a general collection of $n\ge g+1$ points on any ppav $(A, \Theta)$ is theta-general. More
precisely,  a general collection of $n\ge g+1$ points on any non-degenerate  curve $C$ in $A$ is theta-general (see Proposition
\ref{thetasections} for a stronger statement). On the other hand, on a curve $C$ of degree $d=C\cdot\Theta < g$ no collection
of points on $C$ is theta-general: a theta-translate meeting $C$ in $g$ points must contain $C$.

\noindent (b) (Abel-Jacobi curves).  Let $C$ be a curve of genus $g$ and let $A=J(C)$. A {\it general} collection $\Gamma$ of
$n \ge g+1$ points on an Abel-Jacobi image of $C$ is  theta-general. It is interesting to see precisely how this happens. Let $Y\subset C$ be a collection with $|Y|=g$, general in the sense that $h^0(\OO_C(Y))=1$. Although $\Theta^g=g!$, there is a \emph{unique} theta-translate $\Theta_{\gamma_Y}$ such that $\Theta_{\gamma_Y}\cap C=Y$. This is an immediate consequence of the Jacobi inversion theorem. Hence for any other point $p\in C$, $\Theta_{\gamma_Y}$ is the
\emph{unique} theta-translate containing $Y$ and avoiding $p$. This is formalized in the fact that -- denoting $W(Y)$ the locus of theta-translates containing $Y$ (see Notation \ref{notation}) --
$W(Y)$  has two irreducible components, one of which of unexpectedly big dimension. The first component is $W(C)$, the
 locus of all the theta-translates containing the entire curve $C$ (a $(g-2)$-dimensional variety isomorphic to $W_{g-2}$).
  The second one is an isolated point, corresponding to $\Theta_{\gamma_Y}$. This phenomenon actually characterizes Abel-Jacobi curves (see Remark \ref{remark} below)!
\end{example}

\begin{remark}
The notions of theta-independence and theta-generality are considerably weaker than
the corresponding notions  in projective space, essentially due to the large self-intersection of $\Theta$. For example, note that three distinct points on $A$ are always
theta-independent. Moreover, the example above shows that there are theta-general sets
contained in a large family of theta-translates.
\end{remark}

\subsection{Bound on the number of conditions}
The natural analogue of Castelnuovo's basic remark, on the number of conditions imposed on homogeneous forms of given degree,
is given below. Note that for ppav's it is necessary to replace the linear system of hypersurfaces of degree $k$ with the
continuous system formed by all linear systems $|(k\Theta)_\alpha|$, with $\alpha\in\widehat A$.

\begin{proposition}\label{conditions}
 Let $\Gamma\subset A$ be a theta-general set of points and let $k\ge 2$. Then $\Gamma$ imposes at least
$\min\{|\Gamma|,(k-1)g+1\}$ conditions on $H^0(\OO_A((k\Theta)_{\alpha}))$ for $\alpha$ general in $\widehat A$.
\end{proposition}
\begin{proof}
 If $|\Gamma|\ge (k-1)g+1$, the statement means that there $\Gamma$ contains a subset of $(k-1)g+1$
 points imposing independent conditions on $H^0(\OO_A((k\Theta)_{\alpha}))$ for $\alpha$ general in $\widehat A$.
  On the other hand, if $|\Gamma|<(k-1)g+1$, $\Gamma$ can
 be completed to a theta-general subset of $(k-1)g+1$ points. Therefore, it is enough
 to  assume that $|\Gamma|=(k-1)g+1$, and prove that $\Gamma$ imposes independent conditions
  on $H^0(\OO_A((k\Theta)_{\alpha}))$ for $\alpha$ general in $\widehat A$. Let's write
$$\Gamma=Y_1\cup\ldots\cup Y_{k-1}\cup \{p\},$$
with $|Y_i|=g$ for any $i$. \ For any $i$, let  $\gamma_{Y_i,p}$ be
a theta-translate containing $Y_i$ and avoiding $p$. We have that
$\Theta_{\gamma_{Y_1,p}}+\ldots +\Theta_{\gamma_{Y_{k-1},p}}+
\Theta_{\beta}$ contains $\Gamma-\{p\}$ and avoids $p$, unless
$\beta\in \Theta_p$. Therefore,  for any $\alpha\not\in
\Theta_{\gamma_{Y_1,p}+\ldots + \gamma_{Y_{k-1},p}+ p}$, there is a
divisor in $|(k\Theta)_{\alpha}|$ containing $\Gamma-\{p\}$ and
avoiding $p$. As this can be done for any $p$, we have that if
$\alpha\not\in\cup_{p\in \Gamma} \Theta_{\gamma_{Y_1,p}+\ldots +
\gamma_{Y_{k-1},p}+ p}$, then $\Gamma$ imposes independent
conditions on $H^0(\OO_A((k\Theta)_{\alpha}))$.
\end{proof}

\begin{example}[Abel-Jacobi curves]\label{ajexample}
The bound in the previous Proposition is sharp, as seen by looking at Abel-Jacobi curves. Let $C\subset J(C)$ be one such.
Then, for $k\ge 2$ and for any $\alpha\in \widehat{A}$, $h^0(\OO_C((k\Theta)_{\alpha}))= (k-1)g+1.$ Therefore a general
collection of at least $(k-1)g+1$ points on $C$ imposes generically the minimal number of conditions -- namely $(k-1)g+1$ --
on $H^0(\OO_A((k\Theta)_{\alpha}))$. In fact we can be precise: we have the exact sequence
$$ 0\rightarrow H^0(\II_C(k\Theta_{\alpha}))
\rightarrow H^0(\II_{\Gamma}(k\Theta_{\alpha}))\rightarrow H^0(\OO_C(k\Theta_\alpha-\Gamma))\rightarrow 0,$$ where the last
zero follows from the $3$-Theta-regularity of $C$ (cf. \S2(7)). On the other hand, $H^0(\OO_C(k\Theta_{\alpha}-\Gamma))=0$ for $\alpha$ outside
a proper closed subset isomorphic to $W_d$, with $d=kg-|\Gamma|$.
\end{example}

\begin{example}\label{matsusaka}
To illustrate our point of view, let us show how the elementary Proposition \ref{conditions} has as an
immediate consequence the following easy but important step in the proof of the Matsusaka-Ran
criterion: \emph{ Let $C$ be a non-degenerate irreducible curve of degree $g$ in $A$. Then $p_a(C)\le g$.} (It follows
that $C$ is smooth, irreducible, of genus $g$.)

Indeed, let $k$  such that $h^1(\OO_C((k\Theta)_\alpha))$ is generically
zero. Let $\Gamma$ be a general collection of distinct
points on $C$, with $|\Gamma|$ big enough so that $h^0(\OO_C((k\Theta)_\alpha-\Gamma))=0$ (e.g. $|\Gamma|>kg$). Therefore
$h^0(\II_C((k\Theta)_\alpha))=h^0(\II_{\Gamma}((k\Theta)_\alpha))$. By Example \ref{exthetageneral}(a) $\Gamma$ is  theta-general,  hence
by Proposition \ref{conditions} it imposes $\ge (k-1)g+1$ conditions on $H^0(\OO_A((k\Theta)_\alpha))$ for $\alpha$ general in
$\widehat A$. This implies $h^0(\II_C((k\Theta)_\alpha))\le h^0(\OO_A((k\Theta)_\alpha))-(k-1)g-1$, and therefore
$h^0(\OO_C((k\Theta)_\alpha))\ge (k-1)g+1$.
Since $h^1(\OO_C((k\Theta)_\alpha))=0$, by Riemann-Roch we get
$h^0(\OO_C((k\Theta)_\alpha))=kg-p_a(C)+1$. Hence $p_a(C)\le g$.
\end{example}

\subsection{Loci of linear dependence}
In this subsection we study the loci of $\alpha\in\widehat A$ such that a given finite set $\Gamma$ fails to impose
independent conditions on $H^0(\OO_A((k\Theta)_\alpha))$.

\begin{definition/notation}
Let $\Gamma$ be a finite set (or scheme) on $A$. We will consider the cohomological support loci
$$V_{r}(\II_{\Gamma}(k\Theta)) :=\ \{ \alpha\in \widehat A \> ~|~  h^1(\II_{\Gamma}((k\Theta)_\alpha)\ge r \}.$$
For $r=1$, we will simply denote
$$V(\II_{\Gamma}(k\Theta)) :=V_{1}(\II_{\Gamma}(k\Theta)).$$
Since $h^i(\OO_A((k\Theta)_\alpha))=0$ for any $i>0$ and $\alpha\in\widehat A$,  $V(\II_{\Gamma}(k\Theta))$ is the locus of
$\alpha$'s such that $\Gamma$ fails to impose independent conditions on $|(k\Theta)_\alpha|$. For example, by Proposition
\ref{conditions}, for a theta-general collection $\Gamma$ of at most $(k-1)g+1$ points,  $V(\II_{\Gamma}(k\Theta))$ is always a proper subvariety.
\end{definition/notation}

\begin{definition/notation}
Let $\Gamma$ be a collection of points in $A$ and let $p\in \Gamma$. We denote
$$B(\II_{\Gamma-\{p\}}(k\Theta),p) :=\{\>\alpha\in\widehat A\>~|~\hbox{$p$ is in the base locus of $|\II_{\Gamma-\{p\}}((k\Theta)_\alpha)|$}\}.$$

\noindent
We have the basic relation
\begin{equation}\label{dependence}
V(\II_{\Gamma}(k\Theta))=B(\II_{\Gamma-\{p\}}(k\Theta),p)\cup V(\II_{\Gamma-\{p\}}(k\Theta))
\end{equation}
which simply means that if $\Gamma$ fails to impose independent conditions on $|(k\Theta)_\alpha|$, then
either $p$ is in the
base locus of $|\II_{\Gamma-\{p\}}((k\Theta)_\alpha)|$ or $\Gamma-\{p\}$ itself fails to impose independent
 conditions. Note also that, while
$V(\II_{\Gamma}(k\Theta))$ and $V(\II_{\Gamma-\{p\}}(k\Theta))$ are
closed, $B(\II_{\Gamma-\{p\}}(k\Theta),p)$ is only locally closed.
We have also a second basic relation
\begin{equation}\label{basepoints}
V(\II_{\Gamma}(k\Theta))=\bigcup_{p\in \Gamma}B(\II_{\Gamma-\{p\}}(k\Theta),p).
\end{equation}
\end{definition/notation}

The next Lemma describes the intersection of the linear dependence loci of two ``close" collections.

\begin{lemma}\label{elementary}
Let $Z,T$ be finite collections of distinct points of the same cardinality, having all points  but one in common (in
other words $|Z|=|T|=n$ and $|Z\cap T|=n-1$). Then
$$V(\II_Z(k\Theta))\cap V(\II_T(k\Theta))=V_2(\II_{Z\cup T}(k\Theta))\cup V(\II_{Z\cap T}(k\Theta)).$$
\end{lemma}
\begin{proof}
We prove that the left-hand side in contained in the right-hand
side. Denote $p$ (resp. $q$) the only point of $Z$ (resp. of $T$)
which does not belong to $Z\cap T$.
 Assume that both $Z$ and $T$
fail to impose independent conditions on
$H^0(\OO_A((k\Theta)_\alpha))$. If $\alpha\not\in V(\II_{Z\cap T})$
-- i.e. if the $(n-1)$ points of $Z\cap T$ impose independent
conditions on $H^0(\OO_A((k\Theta)_\alpha))$ -- then $p$ is in the
base locus of $|I_{Z\cap T}((k\Theta)_\alpha)|$.  Similarly $q$ is
in the base locus of $|I_{Z\cap
T}((k\Theta)_\alpha)|=|I_{Z}((k\Theta)_\alpha)|$. This means that
$\alpha\in V_2(\II_{Z\cup T}(k\Theta))$. The reverse inclusion is
obvious.
\end{proof}

We conclude with a rough estimate for the dimension of loci of linear dependence
of theta-independent collections.

\begin{lemma}\label{codimension}
Let $Y$ be a theta-independent collection of $n\le g$ points on $A$.
Then $$\dim V(\II_Y(2\Theta))\le g-2.$$
\end{lemma}
\begin{proof}
Fix $p\in Y$ and decompose $Y$ as $Y=T\cup\{p\}$. For any $\gamma\in
\H^{T,p}$ and for any theta-translate $\Theta_\alpha$ such that
$p\not\in\Theta_\alpha$, the divisor $\Theta_\gamma+\Theta_\alpha$
contains $T$ and avoids $p$. This means that, for such $\alpha$'s,
$p$ is not in the base locus of $|\II_T((2\Theta)_\alpha)|$. Thus
$B(\II_{T}(2\Theta),p)$ is contained in $ \bigcap_{\gamma\in
\H^{T,p}}\Theta_{\gamma+p}$, which certainly has codimension $\ge 2$
(note that $\H^{T,p}$ has positive dimension and, since $\Theta$ is
assumed to be irreducible, the intersection of any pair of distinct
theta-translates has codimension $2$). Using (\ref{dependence}),
the assertion follows by induction.
\end{proof}

\begin{example} 
The above estimate is sharp for  a collection of $g$ points lying on an Abel-Jacobi curve.
Indeed, arguing as in  Example \ref{ajexample}, $n\le g+1$ points on
an Abel-Jacobi curve $C$ impose independent conditions on
$H^0(\OO_A(2\Theta))$ away from a locus isomorphic to $W_{n-2}$. It
seems likely that, under the hypothesis of the above Lemma, the more
refined inequality \   $\dim V(\II_Y(2\Theta))\le n-2 $ \  holds.
\end{example}

\section{The Castelnuovo-Schottky Lemma: statement and related results}

\subsection{The statement and some consequences.}
In this section we state the main result of the paper, Theorem \ref{alpha}.
Namely, we will characterize the extremal case of the bound on the number of conditions provided by Proposition \ref{conditions}:

\begin{definition}[Extremal position]
A set of $n\ge g+1$ distinct points $\Gamma \subset A$ is in \emph{extremal position} if it is  theta-general, and if
it imposes precisely $g+1$ conditions on $H^0(\OO((2\Theta)_{\alpha}))$ for general $\alpha\in \widehat{A}$. (Note that
by semicontinuity $\Gamma$ will impose at most $g+1$ conditions on all $2$-theta linear series.)
\end{definition}

We will show that the existence of points in extremal position is intimately related to the existence of trisecants to the Kummer variety associated to $A$, and in fact of a special positive-dimensional family of such. As a consequence, the well-known
Gunning-Welters criterion will imply Theorem \ref{alpha}, which we restate here for
convenience:

\begin{theorem}\label{maintext}
Let $(A,\Theta)$ be an irreducible principally polarized abelian variety of dimension $g$, and let $\Gamma$ be a set of
$n\geq g+2$ points on $A$ in extremal position. Then $(A, \Theta)$ is the canonically polarized Jacobian of a curve $C$, and $\Gamma
$ is contained in a unique Abel-Jacobi embedding $C\subset J(C)$.
\end{theorem}

Admitting Theorem \ref{maintext}, there is a natural way of recovering the curve $C$ from the given data. (In particular this provides another proof of the Torelli theorem.)

\begin{corollary}\label{torelli}
In the setting of Theorem \ref{maintext}, let $U_{\Gamma} \subset \widehat A$ be the open set of $\alpha\in\widehat A$ such that $\Gamma$ imposes exactly $g+1$ conditions on $H^0(\OO_A((2\Theta)_\alpha))$. Then, for any
non-empty open subset $U\subset U_{\Gamma}$, the Abel-Jacobi curve $C$ is the (scheme-theoretic) intersection of all
$Q_\alpha\in|\II_{\Gamma}((2\Theta)_\alpha)|$ with $\alpha\in U$.
\end{corollary}
\begin{proof}
If $\alpha\in U_{\Gamma}$, then $H^0(\II_{\Gamma}((2\Theta)_\alpha))=H^0(\II_C((2\Theta)_\alpha))$.
The statement follows then immediately from
Corollary 4.2 of \cite{pp1}, which says that the sheaf $\II_C(2\Theta)$ is \emph{continuously globally generated}, i.e. for any
open set $U\subset \widehat A$ the evaluation map
$$\bigoplus_{\alpha\in U}H^0(\II_C((2\Theta)_\alpha))\otimes\alpha^\vee\rightarrow \II_C(2\Theta)$$
is surjective.
\end{proof}

As another consequence, we have a similar Schottky-type criterion based on conditions imposed on higher order theta functions.

\begin{corollary}\label{higher}
 Let $(A,\Theta)$ be an irreducible principally polarized abelian variety of dimension $g$, and let $\Gamma$ be a theta-general set on $A$. Assume that there is an integer $k\ge 2$ such that $|\Gamma|>(k-1)g+1$, and that $\Gamma$ imposes exactly $(k-1)g+1$ conditions on
$H^0(\OO(k\Theta_{\alpha}))$ for general $\alpha\in \widehat{A}$. Then $(A, \Theta)$ is the canonically polarized Jacobian of a
curve $C$, and $\Gamma$ is contained in a unique Abel-Jacobi embedding $C\subset J(C)$. If $k\ge 3$, given any $\bar\alpha\in\widehat
A$ such that $\Gamma$ imposes $(k-1)g+1$ conditions on $H^0(\OO_A((k\Theta)_{\bar\alpha}))$, then $C$ is the base locus of
$|I_{\Gamma}(k\Theta_{\bar\alpha}))|$.
\end{corollary}
\begin{proof}
For $k=2$ this is just Theorem \ref{maintext}. If $k>2$, we claim
that the hypothesis implies that any subset $X\subset \Gamma$ such
that $|X|=|\Gamma|-(k-2)g$ is in extremal position. The statement
will then follow from Theorem \ref{maintext}. To prove the claim, we
may assume that $|\Gamma|=(k-1)g+2$, so that $|X|=g+2$. If $X$ is
not in extremal position, i.e. if $X$ imposes independent conditions
on $H^0(\OO((2\Theta)_{\alpha}))$ for general $\alpha\in\widehat A$,
then for any $p\in X$ and for general $\alpha\in \widehat A$, there
is a divisor $D_{\alpha,p}\in |(2\Theta)_\alpha|$ such
$D_{\alpha,p}$ contains $X-\{p\}$ and avoids $p$. We decompose
$\Gamma$ as $\Gamma=(X - \{p\})\cup Y_1\cup\ldots\cup
Y_{k-2}\cup\{p\}$, with $|Y_i|=g$. For any $i$ we choose a
theta-translate $\Theta_{\gamma_i,p}$ containing $Y_i$ and avoiding
$p$. Then
$D_{\alpha,p}+\Theta_{\gamma_1,p}+\ldots+\Theta_{\gamma_{k-2,p}}\in
|(k\Theta)_{\alpha+\sum \gamma_i}|$ is a divisor containing
$\Gamma-\{p\}$ and avoiding $p$. Since $\alpha$ varies in a
Zariski-open set of $\widehat A$, $\alpha+\sum \gamma_{i,p}$ also
does. As this can be done for any $p\in \Gamma$, it follows that
$\Gamma$ imposes independent conditions on
$H^0(\OO((k\Theta)_{\beta}))$, for $\beta$ general in $\widehat A$.
The claim is proved. The last part of the statement follows as in
Corollary \ref{torelli}, since by Corollary 4.2 of \cite{pp1},
$\II_C((k\Theta)_\alpha)$ is globally generated for any $k\ge 3$ and
any $\alpha\in\widehat A$.
\end{proof}

\subsection{Relationship with \cite{grushevsky1}}
After completing a first draft of this manuscript, we were informed by Sam Grushevsky that his paper \cite{grushevsky1} contains a result
proved via the analytic theory of theta functions (and which is used for finding equations for the locus of hyperelliptic Jacobians), whose statement is similar to that of Theorem \ref{maintext}. The statement initially formulated in \cite{grushevsky1} was unfortunately not correct, since the hypothesis was too weak. However recently the author has informed us of a correction, which will appear in an erratum \cite{grushevsky2}. The revised statement is the following:

\begin{theorem}(\cite{grushevsky2})
Let $(A, \Theta)$ be an irreducible ppav of dimension $g$, and $A_0,\ldots,A_{g+1}$ distinct points on $A$. Denote by $K$ the
Kummer map associated to $|2\Theta|$, and suppose that for every $z\in A$, the images $K(A_i + z)$ are linearly dependent.
Assume moreover that there exist some $k$ and $l$ such that for $y:= - \frac{A_k+A_l}{2}$ the linear span of the points
 $K(A_i + y)$ is of dimension precisely $g+1$. Then $A$ is the
Jacobian of some curve $C$, and all the points $A_i$ belong to an Abel-Jacobi embedding of $C$.
\end{theorem}

Half of the hypothesis of both this result and Theorem \ref{maintext}  is the same,  and it would be very interesting to discover a direct relationship between the two complete hypotheses.

\section{Proof of the Castelnuovo-Schottky Lemma}

\subsection{Analysis of loci of linear dependence for points in extremal position.}
>From this point on, unless otherwise stated, $(A,\Theta)$ will be assumed to be an \emph{irreducible} ppav. The following result
is the key property satisfied by sets of points in extremal position in an irreducible ppav.

\begin{lemma/notation}\label{keyresult}
  Let $\Gamma$ be a collection of $n\ge g+2$ points in extremal position on $A$, and let $Y\subset \Gamma$ be any subset consisting of $g$ points. Then:
 \begin{enumerate}
 \item For any $s\in\Gamma-Y$, there is a unique theta-translate $\Theta_{\gamma_Y}$ containing $\Gamma$ and avoiding $s$.
Moreover, this theta-translate works for any $s\in\Gamma-Y$. Hence $\Theta_{\gamma_Y}\cap\Gamma=Y$.
\item For any subset $Z\subset \Gamma$ consisting of $g+1$ points,  the linear dependence locus $V(\II_Z(2\Theta))$ is
a theta-translate,  denoted $\Theta_{\alpha_Z}$.
\item In the setting of (2), if $Z=Y\cup\{p\}$, then $\gamma_Y+p=\alpha_Z.$
\item If $T\subset \Gamma$ is a collection of $g-1$ points and $p,q\in \Gamma-T$, then
$\gamma_{T\cup\{p\}}+q=\gamma_{T\cup\{q\}}+p$.
\end{enumerate}
\end{lemma/notation}
\begin{proof}
Without loss of generality, we can assume that $|\Gamma|=g+2$. We choose two points $p,q\in\Gamma$ and  write $\Gamma=Y\cup
\{p,q\}$, \ $Z :=Y\cup\{p\}$. For any $\gamma\in \H^{Y,q}$ and for any $\beta\in\Theta_p-(\Theta_p \cap\Theta_q)$, the
divisor $\Theta_{\gamma}+\Theta_\beta$ contains $Z$ and misses $q$. Therefore, since we are assuming that for any
$\alpha\in\widehat A$ the set $\Gamma$ fails to impose independent conditions on $|(2\Theta)_\alpha|$, we must have that
$\Theta_{\gamma+p}-(\Theta_{\gamma +p} \cap\Theta_{\gamma +q})$ is contained in  the subvariety $V(\II_Z(2\Theta))$. Since
$\Theta$ assumed to be irreducible, this gives
$$\Theta_{\gamma+p}\subset V(\II_Z(2\Theta))$$ for any $\gamma\in
\H^{Y,q}$. On the other hand, for any $\gamma^\prime\in \H^{Y,p}$ and $\beta\not\in\Theta_p$, the divisor
$\Theta_{\gamma^\prime}+\Theta_\beta$ contains $Y$ and avoids $p$. Hence
$$B(\II_Y(2\Theta),p)\subset \Theta_{\gamma^\prime+p}.$$
In conclusion, for any $p\in T$, it follows from (\ref{dependence}) that
$$\Theta_{\gamma+p}\>\subset \> V(\II_Z(2\Theta))\>\subset\> B(\II_Y(2\Theta),p) \cup V(\II_{Y}(2\Theta))\>\subset \> \Theta_{\gamma^\prime+p}\cup V(\II_{Y}(2\Theta))$$
Since, by Lemma \ref{codimension}, $\dim V(\II_{Y}(2\Theta))\le g-2$, we get that $\gamma=\gamma^\prime$, i.e. that
$\H^{Y,p}$ and $\H^{Y,q}$ consist of the same unique point $\gamma_{Y}$. Moreover $B(\II_Y(2\Theta),p)$ is equal to
$\Theta_{\gamma_Y+p}$ and it is the unique divisor contained in $V(\II_Z(2\Theta))$. Thus the fact that $V(\II_T(2\Theta))=
\Theta_{\gamma_{Y}+p}$  follows from (\ref{basepoints}). We have proved the first three points. Finally, (4) follows immediately from (3).
\end{proof}

\subsection{Existence of trisecants and proof of the Theorem}
Let $k$ be the projective map $A\rightarrow {\bf P}(H^0(\OO_A(2\Theta))^\vee)$. Its image  $k(A)$ is the Kummer variety of
$A$. The relation between points in extremal position and trisecants to $k(A)$ is expressed by the following:

\begin{theorem}\label{trisecants}
Let $\Gamma\subset A$ be a collection of $ g+2$ points in extremal position. Let $p,q,s$ be three points in
 $\Gamma$, and write
$\Gamma=T\cup\{p,q,s\}$ (hence $|T|=g-1$). Then, for any $\alpha\in
\H^{T,s}$, we have that
$$\Theta\cap \Theta_{p-q} \subset  \Theta_{p-s} \cup \Theta_{\alpha+s-q-\gamma_{T\cup\{s\}}}.$$
(See Lemma \ref{keyresult} for the definition of
$\gamma_{T\cup\{s\}}$.)  Equivalently, for any
$\xi\in\frac{1}{2}\H^{T,s}_{p-q-\gamma_{T\cup\{s\}}}$, the points
$k(\xi)$, $k(\xi-(p-q))$ and $k(\xi -( p-s))$ lie on a line.
\end{theorem}
\begin{proof}
The equivalence between the second assertion and the first is well
known (cf. e.g. \cite{mumford} p.80 or \cite{beauville} p.104-105).
For the first assertion, let us fix $p,q\in \Gamma$. Combining Lemma
\ref{elementary}, applied to $T\cup\{p,s\}$ and $T\cup\{q,s\}$, and
Lemma \ref{keyresult}, we have that
\begin{equation}\label{intersection}
\Theta_{\gamma_{T\cup\{s\}}+p}\cap \Theta_{\gamma_{T\cup\{s\}}+q}=V_2(\II_{\Gamma}(2\Theta))\cup V(\II_{T\cup\{s\}}(2\Theta)).
\end{equation}
To begin with, we analyze  the last subvariety appearing on the
right hand side of (\ref{intersection}).
 For any $\alpha\in
\H^{T,s}$ and for any $\beta\not\in\Theta_s$, the divisor
$\Theta_\alpha+\Theta_\beta$ contains $T$ and avoids $s$. Therefore,
as in Lemma \ref{codimension}, we get
$$B(\II_T(2\Theta),s)\subset \bigcap_{\alpha\in \H^{T,s}}\Theta_{\alpha+s}.$$
Hence, by (\ref{dependence}) of \S4,
$$V(\II_{T\cup\{s\}}(2\Theta))\subset \Bigl(\bigcap_{\alpha\in \H^{T,s}}\Theta_{\alpha+s}\Bigr)\cup
V(\II_{T}(2\Theta)).$$ Clearly $V(\II_{T}(2\Theta))\subset
V(\II_{\Gamma-\{s\}}(2\Theta))$ Therefore
\begin{equation}\label{Gamma-p-q}
V(\II_{T\cup\{s\}}(2\Theta))\subset \Bigl(\bigcap_{\alpha\in
\H^{T,s}}\Theta_{\alpha+s}\Bigr)\cup V(\II_{\Gamma-\{s\}}(2\Theta)).
\end{equation}
Now we turn our attention to the first subvariety of the right hand
side of (\ref{intersection}). For any point $r\in \Gamma$ we have
$V_2(\II_{\Gamma}(2\Theta))\subset V(\II_{\Gamma-\{r\}}(2\Theta))$.
In particular
$$V_2(\II_{\Gamma}(2\Theta))\subset
V(\II_{\Gamma-\{s\}}(2\Theta))$$ Putting together with
(\ref{Gamma-p-q}) it turns out that, for any $\alpha\in \H^{T,s}$,
\begin{equation}\label{partial-decomposition}\Theta_{\gamma_{T\cup\{s\}}+p}\cap
\Theta_{\gamma_{T\cup\{s\}}+q}\subset
V(\II_{\Gamma-\{s\}}(2\Theta))\cup \Theta_{\alpha+s}.
\end{equation}
To conclude the proof, note that, for any $r\in \Gamma$,
$V(\II_{\Gamma-\{r\}}(2\Theta))=\Theta_{\alpha_{\Gamma-\{r\}}}$ for
any $r\in\Gamma$ (Lemma \ref{keyresult}).
 If
$p\ne r$, we can write
$$\alpha_{\Gamma-\{r\}}=\gamma_{T\cup\{s\}}+p+q-r,$$
since by Lemma \ref{keyresult} we have that
$\alpha_{\Gamma-\{r\}}=\gamma_{\Gamma-\{p,r\}}+p$ and
$\gamma_{\Gamma-\{p,r\}}+r=\gamma_{\Gamma-\{p,q\}}+q=\alpha_{\Gamma-\{p\}}$
(note that $\Gamma-\{p,q\}=T\cup\{s\}$). In conclusion, for any
$r\in \Gamma$,
\begin{equation}\label{X}
V(\II_{\Gamma-\{r\}}(2\Theta))= \Theta_{\gamma_{T\cup\{s\}}+p+q-r}.
\end{equation}
Plugging this for $r=s$ into (\ref{partial-decomposition}) we get
that, for any $\alpha\in \H^{T,s}$
$$\Theta_{\gamma_{T\cup\{s\}}+p}\cap \Theta_{\gamma_{T\cup\{s\}}+q}\subset\Theta_{\gamma_{T\cup\{s\}}+p+q-s}\cup
\Theta_{\alpha+s}.$$
The statement follows by translating by $-\gamma_{T\cup\{s\}}-q$.
\end{proof}

Now can finally put everything together in order to prove our Castelnuovo-Schottky Lemma for abelian varieties.

\begin{proof} (of Theorem \ref{alpha}).
We will use the following part of Welters' criterion (building on previous work of Gunning \cite{gunning}):

\noindent \emph{Let a,b,c be three distinct points on an irreducible
ppav $(A,\Theta)$. If the locus $W_{a,b,c}$ of $\xi\in A$ such that
$k(2\xi+a)$, $k(2\xi+b)$ and $k(2\xi+c)$ lie on a line in ${\bf
P}(H^0(\OO_A(2\Theta))^\vee)$ is positive-dimensional, then
$W_{a,b,c}$ is a smooth irreducible curve and $(A,\Theta)$ is the
Jacobian of $W_{a,b,c}$.} (\cite{welters2}, Theorem (0.5), case (i);
see also \cite{beauville}, p.104-105).

To prove Theorem \ref{alpha} we fix three points $p,q,s\in \Gamma$
and we consider any subset $X$ of $ \Gamma$, with $|X|=g+2$, and
containing $p,q,s$. We write $X=T\cup\{p,q,s\}$. Since $|T|=g-1$,
(every component of) $\H^{T,s}$ is positive-dimensional (see
Definition/Notation \ref{independence}). Therefore, by Theorem
\ref{trisecants} and the Gunning-Welters criterion, we are on a
Jacobian.  It remains to prove that $\Gamma$ is contained in an
Abel-Jacobi curve. Note that it turns out that the closure
$(\overline{\H^{T,s}})_{p-q-\gamma_{T\cup\{s\}}}$ coincides with
$W_{0,p-q,p-s}$,  hence it is a smooth irreducible curve and
$(A,\Theta)$ is its polarized Jacobian. Moreover, it is easy to
deduce (cf. e.g. \cite{beauville}, p.104-105) that
$(\overline{\H^{T,s}})_{p-q-\gamma_{T\cup\{s\}}}=C_{p-q-s}$, for a
fixed Abel-Jacobi curve $C$ in $A=J(C)$. In particular, it follows
that $$(\overline{\H^{Y,s}})_{\gamma_{T\cup\{s\}}-s}=C$$ does not
depend on $s$ and $T$ but just on $p$ and $q$. By Lemma
\ref{keyresult}(4), $\gamma_{T\cup\{s\}}-s=\gamma_{T\cup\{t\}}-t$,
so
$$\overline{\H^{T,s}}=C_{\gamma_{T\cup\{t\}}-t}.$$
Now if $t\in \Gamma$, $t\ne s$ then $\gamma_{T\cup\{t\}}\in\H^{T,s}$ (Remark
\ref{k-1}). Hence $t\in C$. Therefore $\Gamma\subset C$.
\end{proof}

\section{Genus bound}

As another application of this point of view, we prove Theorem \ref{genus_bound}. This is
a ``Castelnuovo bound", i.e. a bound on the genus of a curve on a $g$-dimensional ppav $(A, \Theta)$ in function of its degree $d: = C\cdot \Theta$.  The
bound is quadratic in the degree, with leading term $d^2/2g$. Although the proof shows,
somewhat subtly, that it is not optimal for $g\ge 3$,
it is of the expected order of magnitude, and it improves considerably a previously known bound  \footnote{In \cite{debarre1}
Theorem 5.1, it is shown that the Castelnuovo bound for curves in projective space yields a quadratic bound for ppav's whose leading term is $2d^2/(g-1)$.}.

\begin{remark}\label{effectivity}
On abelian surfaces, i.e. for $g=2$, the bound in Theorem \ref{genus_bound} is optimal since for even $d =2m$ it reads
$\gamma\le d^2/4+1$, which is the genus of smooth curves in $|m\Theta|$.
\end{remark}

\subsection{General position of general theta-sections of an irreducible curve.}
To prove Theorem \ref{genus_bound}, we will
follow Castelnuovo's method based on the number of conditions imposed on hypersurfaces by a general hyperplane section of the curve $C$ (\cite{acgh}, p.114--115). The main point here is that, as soon as the degree is higher than $g$, a general theta-section
$C\cap\Theta_\alpha$ is theta-general (Proposition \ref{thetasections}).  Although we are able to prove this by means of a direct
argument, we propose the following problem, which is interesting on its own sake.

\begin{problem}
Study the monodromy of the general theta-section of  non-degenerate curves in  irreducible ppav's. (On \emph{general} ppav's it can be shown \` a la Harris --
cf. \cite{acgh} p.111 --  that the monodromy
of a general theta-section is always the symmetric group).
\end{problem}

Here we will confine ourselves to a more elementary analysis, which will be enough for our purpose, and won't tackle the problem above.  Let $C$ be
a non-degenerate, reduced and irreducible curve on $A$. Then it is well-known that there exists a non-empty open set $U\subset
\widehat A$ such that  for any $\alpha\in U$ the theta-translate $\Theta_\alpha$ meets $C$ transversally. We make a
preliminary observation (borrowed from \cite{gh}, Lemma at p.249). Let $U$ be the above open set. Denoting by $C^{(d)}$ the symmetric product of $C$, we have the map
$$\phi:U\rightarrow C^{(d)},~~ \beta\mapsto \Theta_\beta\cap C,$$
whose image does not meet the diagonals. Since $C$ is
non-degenerate, $\phi$ is finite. For a fixed $\alpha\in U$, we put
an order on $\Theta_{\alpha}\cap C=\{p_1,\dots ,p_d\}$. Up to
restricting $U$, we can lift $\phi$ to a map to the cartesian
product $\psi: U\rightarrow C^d$ such that, for $\beta\in U$,
$\psi(\beta)=\{(\Theta_\beta\cap C)_1,\dots , (\Theta_\beta\cap
C)_d\} $. Fix $k\le g$. For any multi-index $I=\{i_1,\dots
,i_{k}\}$, by composing with the corresponding projection, we get a
map $\pi_I:V\rightarrow C^{k}$. The following Lemma follows
immediately.

\begin{lemma}\label{g-uniformity}
The maps $\pi_I$ are dominant for any $I$. Therefore, if $\alpha$ is sufficiently general,
any property satisfied by a general  effective divisor of degree $k\le g$ on $C$
is satisfied by \emph{any} effective divisor $Y$ of degree $k$ contained in $\Theta_\alpha\cap C$.
\end{lemma}

So far for uniformity properties of a general theta-section. Concerning theta-general position, we start with the following:

\begin{lemma}\label{twothetas}
Let $C$ be a non-degenerate, reduced and irreducible curve of degree $d>g$.
Then a general divisor of degree $g$ on $C$ is contained in at least two distinct theta-sections.
\end{lemma}
\begin{proof}
To simplify the notation, we will prove the result assuming that $C$ is smooth (the same argument works in the
 non-smooth case via passage to the
normalization of $C$). Assume that the assertion is not true. Then one can associate to a general
 divisor of degree $g$, say $Y$,
the linear equivalence class on $C$ of the unique theta-section containing $Y$. This induces a rational map
$f:C^{(g)}\rightarrow {\rm Pic}^0(C)$.  Now $f$ has to factor through the Albanese map of $C^{(g)}$, i.e. the Abel-Jacobi map
$C^{(g)}\rightarrow J(C)$. Let $h$ be the induced (endo)morphism of abelian varieties $h:J(C)\rightarrow {\rm Pic}^0(C)$. Note
that, by construction, all the $Y$'s contained in a given theta-section $\Theta_a\cap C$ are contained in a fiber of $f$.
Therefore, as we are assuming $C\cdot\Theta>g$, we can choose two distinct divisors $Y_1$ and $Y_2$ both contained in
$\Theta_a\cap C$, of the form $Y_1=p_1+\ldots +p_{g-1}+p$ and $Y_2=p_1+\ldots +p_{g-1}+q$.  We have that
$h(p)-h(q)=f(Y_1)-f(Y_2)=0$. Since this can be done for general $p$ and $q$ in $C$, we would have that $C$ is contracted by
$h$, i.e. that $h$ is constant. This yields that $f$ is constant, a contradiction.
\end{proof}
\begin{remark}\label{remark} The above Lemma provides another characterization of Jacobians and Abel-Jacobi curves:
\emph{if there exists a non-degenerate curve $C$ such that, given a collection $Y$ of $g$ general distinct points on $C$,
there is only one theta-translate containing $Y$ and not containing $C$, then the abelian variety $A$ is a Jacobian and $C$
is an Abel-Jacobi curve.} This follows at once from the previous Lemma and the Matsusaka-Ran criterion.
\end{remark}

We are now in a position to prove the main technical result of this subsection.

\begin{proposition}[General position]\label{thetasections}
Let $C$ be a non-degenerate, reduced and irreducible curve of degree $d >g$
and let $X=\Theta_\alpha\cap C$ be a general theta-section of $C$. Then $X$ is theta-general.
\end{proposition}
\begin{proof}
Given a divisor $Y\subset X$ such that ${\rm deg}(Y) = g$, we denote by $W(Y)$ the locus of theta-translates
containing $Y$. Moreover we denote $\phi(Y) := C\cap\bigcap_{\alpha\in W(Y)}\Theta_\beta$. By Proposition \ref{g-uniformity}, the
cardinality of $\phi(Y)$ is constant for all such $Y$, and we will denote it by $n$. By definition, the fact that $X$ is theta-general means that $Y=\phi (Y)$ for all $Y$, i.e. that $n=g$.

We make the following claim: \emph{let $Z$ be another divisor of degree $g$ contained in $X$. Then $\phi(Y)= \phi(Z)$ if and only if $Z\subset
\phi(Y)$.} If the converse implication were not true, then $\deg\phi(Z)<\deg\phi(Y)$, contradicting Proposition \ref{g-uniformity}. The direct implication follows from the definition.

Let's denote by ${\mathcal P}^j(X)$ the set of subsets of $X$ of cardinality $j$.
We have produced a family $\Phi=\bigl\{\phi(Y)\bigr\}_{Y\in{\mathcal P}^g( X)}
\subset{\mathcal P}^n(X)$ with the following property: \emph{for any
 $Y\in{\mathcal P}^g(X)$ there exists a \emph{unique} $\phi\in\Phi$ containing $Y$}.
 It follows easily that
 $\Phi$ falls into one
 of the following three cases: \ (1) $n=g$ and $\Phi={\mathcal P}^g(X)$; \ (2) $g=1$, $|X|$ is a multiple of $n$, and $\Phi$ is a partition of
 $X$ in subsets of cardinality $n$; \ (3) $n=|X|$ and $\Phi=\{X\}$
\footnote{Proof. We prove that if $n> g$ then either (2) or (3) hold. If $n=g+1$,
 for any $Y\in{\mathcal P}^g(X)$ the residual set $r_Y=\phi(Y)-Y$ is a point. This
 establishes a map $r:{\mathcal P}(X)\rightarrow X$, $Y\mapsto r_Y$. By construction, the
 map $r$ is bijective. Therefore $|{\mathcal P}^g(X)|=|X|$ and the assertion follows.
 If $n>g+1$, subtracting one point to any $\phi\in\Phi$, we get a new set $X^\prime$,
 equipped with a family of subsets
 $\Phi^\prime\subset {\mathcal P}^{n-1}(X^\prime)$ with the same property. Therefore the assertion
 follows easily by induction on $n$.}. But case (2) is
 excluded since  $g\ge 2$, and case (3) is excluded since, by Lemma \ref{twothetas},
 $n<|X|$. Therefore case (1) holds, and the Proposition is proved.
 \end{proof}

\subsection{Proof of the bound. }
We are now ready for the proof of the genus bound. From this stage on, the argument (for the first part) is essentially that of Castelnuovo, as accounted e.g. in \cite{acgh}, p.115.

\begin{proof} (of Theorem \ref{genus_bound}).
Let $\Theta_{\bar \alpha}$ be a fixed general theta-translate, so that $X=C\cap \Theta_{\bar\alpha}$ is
 theta-general. Let us
denote by $\beta_l$ the generic value of the (affine) dimension of the linear series cut out by $|(l\Theta)_\alpha|$ on $C$,
with $\alpha\in \widehat{A}$. In other words, $\beta_l$ is the generic value of the difference
$$h^0(\OO_A((l\Theta)_\alpha))-h^0(\II_{C}((l\Theta)_\alpha)).$$
Taking $H^0$'s in the exact sequence
\begin{equation}\label{referee}
0\rightarrow \II_C(((l-1)\Theta)_{\alpha-\bar\alpha})\rightarrow
\II_C((l\Theta)_{\alpha})\rightarrow
\II_{X/\Theta_{\bar\alpha}}((l\Theta)_{\alpha})\rightarrow 0
\end{equation}
shows, after an immediate computation,  that for general $\alpha$
and for any $l\ge 1$, the difference $ \beta_l-\beta_{l-1}$ is
greater than or equal to the number of conditions imposed by $X$ on
$H^0(\Theta_{\bar\alpha},\OO_{\Theta_{\bar\alpha}}(l\Theta_\alpha))$,
which is in turn equal to the number of conditions imposed by $X$ on
$H^0(\OO_A(l\Theta_{\alpha}))$. Assume $d>g$. By Proposition
\ref{thetasections}, $X$ is theta-general. Hence, by Proposition
\ref{conditions},
\begin{equation}\label{difference}
\beta_l-\beta_{l-1}\ge \min\{d,(l-1)g+1\}.
\end{equation}
Let $\pi:C^\prime\rightarrow C$ be the birational morphism in the statement, and let us denote by
 $\lambda_h$ the generic value
of $h^0(C^\prime,\pi^*\OO_{C}(h\Theta_\alpha))$. It is clear that
$\lambda_h\ge \beta_{h}$ for any $h$. Let $m
:=\bigr[\frac{d-1}{g}\bigl]$, so that $d-1=gm+\epsilon$, with
$0\le\epsilon<g$. By (\ref{difference}) it follows (noting that
$\beta_1=1$) that for any $k\ge 0$:
\begin{equation}\label{betam}
\beta_{m+1+k}-1=\sum_{l=2}^{m+1+k}(\beta_l-\beta_{l-1})\ge\sum_{l=1}^{m}(lg+1)+kd={{m+1}\choose 2}g+m+kd.
\end{equation}
But, for $k$ sufficiently big,
$h^1(C^\prime,\pi^*\OO_{C}((m+1+k)\Theta_\alpha))$ vanishes.
Therefore, by Riemann-Roch on
 $C^\prime$:
\begin{equation}\label{inequality}
(m+1)d-\gamma+1+kd=\lambda_{m+1+k}\ge \beta_{m+1+k}\ge{{m+1}\choose 2}g+m+1+kd.
\end{equation}
The inequality of the statement follows.

To prove the last part note that, if equality  is attained by a
certain curve $C$ in $A$, then for any $l\ge 2$ we must have (using
(\ref{referee}) again)
\begin{equation}\label{equality}
\beta_l-\beta_{l-1}=h^0(\OO_A((l\Theta)_\alpha)-h^0(\II_X((l\Theta)_\alpha))={\rm min}\{d,(l-1)g+1\}
\end{equation}
for $\alpha$ general in $\widehat A$. From the first equality in (\ref{equality}) it follows easily that if a curve $C$ attains equality, then
\begin{equation}\label{vanishing}
h^1(\II_C((l-1)\Theta_\alpha))=0 \end{equation}  for any $l\ge 2$
and for general $\alpha\in\widehat A$. Given a Zariski-open set
$U\subset \widehat A$, we denote by $S_{l,U}$ (resp. $\Gamma_{l,U})$
the intersection of all $Q_\alpha\in
|\II_{C}(((l+1)\Theta)_{\alpha+\bar\alpha})|$ (resp. the
intersection of all $E_\alpha\in
|\II_{X}(((l+1)\Theta)_{\alpha+\bar\alpha})|$), for $\alpha\in U$.
>From the exact sequence (\ref{referee}) and from (\ref{vanishing})
it follows that, for a suitable $U\subset \widehat A$,
$\Gamma_{l,U}=S_{l,U}\cap \Theta_{\bar\alpha}$. But if $d\ge
(l-1)g+2$, Theorem \ref{maintext} and Corollaries \ref{torelli} and
\ref{higher}, together with the second equality in (\ref{equality}),
yield that $A$ is a Jacobian and $\Gamma_{l,U}$ is an Abel-Jacobi
curve. This would imply that $S_{l.U}$ is a surface whose generic
theta-section is an Abel-Jacobi curve, which is impossible if $g\ge
3$. In conclusion, (\ref{vanishing}) must fail for any $l\ge 2$ such
that $(l-1)g+1<d$. Therefore (\ref{equality}) must fail for at least
one $l$, so there is strict inequality in the genus bound.
\end{proof}

\section{Appendix: the divisor class associated to  points failing to impose independent conditions}

In Proposition \ref{keyresult} it is shown, by elementary methods, that given a collection $\Gamma$ of $g+2$ or more points in \emph{extremal position} then, for any subset $Z\subset\Gamma$ with $|Z|=g+1$, the locus $V(\II_Z(2\Theta))$ is a
\emph{divisor} (in fact a specific theta-translate). This is a key point in the proof of
Theorem \ref{alpha}. It turns out that a weaker version of this statement holds in great generality. The proof uses the $M$-regularity criterion (cf. \cite{pp1}).

Let $A$ be an abelian variety, $D$ an ample divisor, and $Z$ a finite set on $A$.  Recall that we denote
$$V(\II_Z(D)) : = \{\alpha\in\widehat A\>|\>h^1(\II_Z(D_\alpha))>0\}.$$

\begin{proposition}\label{divisor}
Let $\Gamma$ be a collection of distinct points on $A$. Let $n(\Gamma)$ be  the number of conditions generically imposed by
$\Gamma$ on $H^0(\OO_A(D_\alpha))$. Assume that $n(\Gamma)<|\Gamma|$  and let $Z\subset \Gamma$ be a subset of $n(\Gamma)$
points generically imposing independent conditions on $H^0(\OO_A(D_\alpha))$. Then $V(\II_Z(D))$ is a proper closed subset of
$\widehat A$, containing at least one divisorial component.
\end{proposition}
\begin{proof}
The fact that $V(\II_Z(D))$ is a proper subvariety of $\widehat A$ follows by definition. To prove that it contains a
codimension one component, we note that since $Z$ is finite and $D$ is ample, $h^i(\II_Z(D_\alpha))=0$ for all $i\ge 2$ and
all $\alpha\in \widehat A$. Therefore $\II_Z(D)$ is an \emph{$M$-regular sheaf} if and only if
${\rm codim}_{\widehat A}V(\II_Z(D))\ge
2$ (cf. \cite{pp1} \S2). If this is the case,
 the $M$-regularity criterion (Corollary 4.2 of \emph{loc. cit.}) yields that $\II_Z(D)$ is
\emph{continuously globally generated}. This means that, for any open set $U\subset \widehat A$ the evaluation map
$$\bigoplus_{\alpha\in
U}H^0(\II_Z(D_\alpha))\otimes\alpha^\vee\rightarrow H^0(\II_Z(D))$$ is surjective. It follows immediately that, if ${\rm
codim}_{\widehat A}V(\II_Z(D))\ge 2$, then for $\alpha$ general in $\widehat A$ there is no subscheme $\Gamma$ containing $Z$
strictly and such that $H^0(\II_Z(D_\alpha))=H^0(\II_\Gamma(D_\alpha))$.
\end{proof}

Assume that a subset $\Gamma$ as in the previous Proposition is in a sufficiently \emph{uniform position}  (we won't give a precise
definition here). Then, for any subset $Z\subset \Gamma$ with $|Z|=n(\Gamma)$, the subvarieties $V(\II_Z(D_\alpha))$ are
proper and the algebraic equivalence classes of their divisorial part coincide. It is possible that the divisorial parts
of the $V(\II_Z(D_\alpha))$'s have common components. We call the class the remaining components the (mobile) \emph{divisorial class} of $\Gamma$. If $\Gamma$ is a general subset of points on a curve $C$, the divisorial class of $\Gamma$ is related to $C$.

\begin{example}\label{example}
Let $(A,\Xi)$ be ppav and let $C$ curve in $A$ such that: \ (a)
 $h^1(\II_C((l\Xi)_\alpha))=0$ for $\alpha$ general in $\widehat A$, \
  (b)
 $h^1(\OO_C(l\Xi)_\alpha))=0$ outside a subvariety of codimension
 at least two in $\widehat A$.  \
 In the previous notation, take $D=l\Xi$. Let $\Gamma$ be a general effective
divisor on $C$, of sufficiently high degree, so that $h^0(\OO_C((l\Xi)_\alpha-\Gamma))=0$. Then, using (a), $n(\Gamma)=
 h^0(\OO_C(l\Xi))=l\cdot d(C)-g(C)+1$. Let $Z\subset \Gamma$ be a general  effective divisor of degree $n(\Gamma)$ and
 let $\Theta_{Z}$ be the theta-translate (in the Jacobian of $C$) given by
$\{\beta\in {\rm Pic}^0(C)\>|\>h^1(\OO_C((l\Xi)_\beta-Z)>0\}$ (we have $l\Xi\cdot C-n(\Gamma)=g(C)-1$). We consider
  the map $\pi: \widehat A\rightarrow {\rm Pic}^0(C)$. We have the exact sequence
 \begin{equation}\label{last}
 H^1(\II_C((l\Xi)_\alpha))\rightarrow  H^1(\II_{Z}((l\Xi)_\alpha))\rightarrow
 H^1(\OO_C((l\Xi)_\alpha-Z))\rightarrow  H^2(\II_C((l\Xi)_\alpha))
 \end{equation}
 Since $H^2(\II_C((l\Xi)_\alpha))\cong H^1(\OO_C(l\Xi)_\alpha))$,  (b) ensures that, regarding divisorial
 components of  $V(\II_Z(l\Xi))$, the last term of  (\ref{last}) is neglectable.
 Moreover (\ref{last}) yields that the divisorial part of
 $V(\II_Z(l\Xi))$ is $\pi^* \Theta_Z$, plus a (possibly empty) fixed divisor contained
 in $V(\II_C(l\Xi))=\{\alpha\in\widehat A\>|\> h^1(\II_C((l\Theta)_\alpha)>0\}$. Thus the (mobile) divisorial class of
 $\Gamma$ is $\pi^*[\Theta]$. For example, let $\tilde C$ be an Abel-Prym curve in a Prym variety, and $l=2$.
 As  E. Izadi informs us, it can be deduced from \cite{izadi} that $h^1(\II_{\tilde C}((2\Xi)_\alpha)$ vanishes
  for general $\alpha\in
\widehat A$. Hence (a) holds, while (b) holds trivially. Therefore the divisorial class of a general
such $\Gamma\subset \tilde C$ is $2[\Xi]$.
\end{example}

We believe that, in view of possible extensions of the
Castelnuovo-Schottky Lemma as in \S2, the consideration of the divisor class is crucial. To start with, given an irreducible
ppav $(A,\Theta)$, it would be interesting to know a lower bound for the number $n(\Gamma)$ of a uniform, theta-general
collection $\Gamma$ with associated class $n[\Theta]$ (cf. the end of \S2 above for a conjecture in the case
 $n=2$).

\providecommand{\bysame}{\leavevmode\hbox to3em{\hrulefill}\thinspace}

\end{document}